\documentclass[oneside,11pt]{article} 

\usepackage{epsfig,graphicx}
\usepackage{latexsym,amssymb}
\usepackage{setspace,cite} 

\usepackage{amsmath, amssymb, amsthm}
\usepackage{graphicx,color}

\usepackage{fancyhdr}

\usepackage{a4}
\usepackage{bbm}

\newcommand{\B}[1]{\mathbb{#1}}
\usepackage{amsmath}
\usepackage{amssymb}
\usepackage{euscript}
\usepackage{eufrak}
\usepackage{amsthm}
\usepackage[all]{xy}
\newtheorem{prop}{Proposition}[section]

\newtheorem{cor}[prop]{Corollary}
\newtheorem{lemma}[prop]{Lemma}
\newtheorem{thm}[prop]{Theorem}
\newtheorem{conj}[prop]{Conjecture}

\theoremstyle{definition}
\newtheorem{defn}[prop]{Definition}

\newtheorem{queste}[prop]{Question}
\theoremstyle{definition}
\newtheorem*{rmk}{Remark}
\theoremstyle{definition}

\theoremstyle{definition}

\title{Higher Derivatives of L-series associated to Real Quadratic Fields}
\author{Lawrence Taylor}

\begin{document}
\maketitle

\begin{abstract}
This text is a modified version of a chapter in a PhD thesis \cite{TaylorPhD} submitted to Nottingham University in September 2006, which studied an approach to Hilbert's twelfth problem inspired by Manin's proposed theory of Real Multiplication \cite{Manin}. In \cite{TaylorII} we defined and studied a nontrivial notion of line bundles over Quantum Tori.  In this text we study sections of these line bundles leading to a study concerning theta functions for Quantum Tori.  We prove the existence of such meromorphic theta functions, and view their application in the context of Stark's conjectures and Hilbert's twelfth problem.  Generalising the work of Shintani, we show that (modulo a Conjecture \ref{conjh}) we can write the derivatives of L-series associated to Real Quadratic Fields in terms of special values of theta functions over Quantum Tori.
\end{abstract}

\section{Introduction}

The theory of Complex Multiplication uses the theory of elliptic curves to provide a complete explicit class field theory for imaginary quadratic fields, thereby solving Hilbert's twelfth problem for such fields.  The proposed theory of Real Multiplication \cite{Manin} is expected to solve the corresponding problem for real quadratic fields.\newline

As discussed in \cite{Taylor,TaylorII}, non-Hausdorff topological spaces known as Quantum Tori are supposed to play an analogous role in Real Multiplication as elliptic curves do in Complex Multiplication.  In \cite{TaylorII} we developed a nontrivial notion of line bundles over these spaces, and in this text we discuss the existence of sections of these nontrivial bundles.\newline

Sections of line bundles over topological spaces are a source of many interesting functions in number theory. Examples include the modular $j$-function, modular forms, and Jacobi theta functions.  Given the expected role of Quantum Tori in an explicit class field theory for real quadratic fields, we expect sections of line bundles over these objects to be interesting and relevant to such a theory.  We will see that this is indeed so through the work of Shintani, and together with Stark's conjectures may provide crucial information on generating abelian extensions of real quadratic fields.\\

In \cite{TaylorII} we defined the notion of a line bundle over a Quantum Torus $Z_L$ as an element of the group of cocycles $Z^1(L,\mathcal{H}^\ast)$.   The motivation for this definition was the fact that line bundles over a Complex Torus $X_\Lambda$ are in bijection with the group of cocycles $Z^1(\Lambda,\mathcal{H}^\ast)$.  Suppose $\pi_\mathcal{L}:\mathcal{L} \rightarrow X_\Lambda$ is a line bundle over such a Complex Torus $X_\Lambda$, and corresponds to a cocycle $A_\lambda(v) \in Z^1(\Lambda,\mathcal{H}^\ast)$.  The topological space $\mathcal{L}$ is viewed as the quotient of $\B{C} \times \B{C}$ by the action of $\Lambda$ given by 
$$\lambda(z,v)=(z+\lambda,A_\lambda(v)z).$$
A section of $\mathcal{L}$ is a map $\sigma : X \rightarrow \mathcal{L}$, such that $\pi_\mathcal{L} \circ \sigma=1_{X_\Lambda}$. If $p$ and $\tilde p$ denote the natural projections
$$
\xymatrix{
\B{C} \ar[r]& X_\Lambda & \qquad \textrm{and}\\
\B{C} \times \B{C} \ar[r] & \mathcal{L} & 
}
$$
respectively, then we have the following commutative diagram:
$$
\xymatrixcolsep{5pc}
\xymatrix{
\B{C} \times \B{C} \ar[d]^\pi \ar[r]^{\tilde p}& \mathcal{L} \ar[d]_{\pi_{\mathcal{L}}}\\
\B{C} \ar[r]^p& X_\Lambda \ar@/_1pc/[u]_\sigma.
}
$$
where $\pi$ is the projection on to the first coordinate.\\

The natural projection $\tilde p:\B{C} \times \B{C} \rightarrow \mathcal{L}$
is a covering map, so the section $\sigma$ lifts to a section $\tilde \sigma: \B{C} \rightarrow \B{C} \times \B{C}$ of the trivial bundle on $\B{C}$ satisfying 
\begin{equation}\label{pisigma}
\pi \circ \sigma=1_{\B{C}}.
\end{equation} 
By \eqref{pisigma} we have
$$\tilde \sigma(z)=(z,\theta(z))$$
for some $\theta \in \mathcal{H}^\ast$.  Since $\tilde \sigma$ is a lift of $\sigma$, for all $\lambda \in \Lambda$ we have
$$ \tilde p \circ \tilde \sigma(z)=\tilde p \circ \tilde \sigma(z+\lambda),$$
which implies that
\begin{equation}\label{periodtheta}
\theta(z+\lambda)=A_\lambda(v)\theta(z).
\end{equation}
Conversely, if $\theta \in \mathcal{H}^\ast$ satisfies the periodicity condition of \eqref{periodtheta} with respect to the lattice $\Lambda$, the map 
$$\sigma: z+\Lambda \mapsto \tilde p((z,\theta(z)))$$
defines a section of $\mathcal{L}$.  This prompts the following definition:

\begin{defn}[Theta function]
A holomorphic theta function for a group $G \subseteq \B{C}$ is a holomorphic map $\theta: \B{C} \rightarrow \B{C}$
such that for all $v \in \B{C}$
\begin{equation} \label{thetacondition}
\theta(v+g)=A_g(v)\theta(v)
\end{equation} for some $A_g \in Z^1(G,\mathcal{H}^\ast)$.
\end{defn}

Hence theta functions for a complex lattice $\Lambda$ correspond to sections of holomorphic line bundles over the Complex Torus $X_\Lambda$.  To determine the existence of sections of holomorphic line bundles over Quantum Tori $Z_L$, we need to determine whether there
are any holomorphic theta functions for the pseudolattice $L$.  \\

This text is split in to two main parts, the first consisting of \S\ref{holo} and \S\ref{mero}.  In the first of these we show that unlike the case for complex lattices $\Lambda$, there are no nontrivial holomorphic theta functions corresponding to a pseudolattice $L$.  In \S\ref{mero}, we weaken the condition of holomorphicity to allow our theta function to have poles. We show that the double sine function studied by Shintani \cite{ShintaniII,ShintaniIII} and Kurokawa \cite{KurokawaII,KurokawaI} can be interpreted to be a \emph{meromorphic} theta function for a pseudolattice.\\

In the 1970's Stark made a series of conjectures \cite{StarkI,StarkII,StarkIII,StarkIV} regarding the values of L-functions associated to number fields at $s=0$.  The second half of this text concerns the application of the functions we discuss in the first part to these conjectures.  In \S\ref{starks conj} we give an introduction to Stark's ideas, and how they are related to our goal of understanding an explicit class field theory for real quadratic fields.  The remainder of \S\ref{starks conj} is devoted to an account of the work of Shintani.  In \cite{ShintaniII} Shintani described the values of an L-function associated to a real quadratic field in terms of specific values of the double sine function, and in a later paper \cite{ShintaniIII} used these values to generate abelian extensions of specific real quadratic fields.  In the context of this thesis, this is an important result, stating that in specific cases special values of \emph{meromorphic} theta functions associated to Quantum Tori can generate abelian extensions of certain real quadratic fields.\\

We see that Shinatani's result can be interpreted as a solution in a special case to the Rank One Abelian Stark conjecture ~\cite{Roblot,Tangedal}, which concerns the case when the L-function has a simple zero at $s=0$.  There exist higher order Stark conjectures ~\cite{Tate,Rubin} which concern the cases when the L-function has zeros of higher order at $s=0$.  Motivated by these conjectures and Shintani's result we investigate whether it is possible to write higher derivatives of L-functions associated to real quadratic fields in terms of meromorphic theta functions for a pseudolattice. Our main result is Theorem \ref{special}, which writes the $m^{th}$ derivative of an 
L-function as an element of a certain field, whose generators contain the special values of various functions defined in \S\ref{gam} which are shown to be theta functions for pseudolattices. \\

We use Shintani's work of \cite{ShintaniII} to reduce the proof of Theorem \ref{special} to a result concerning a type of zeta function.  This result is proved in \S\ref{proof} using a blend of induction (of which Shintani's result is the starting case), and the calculation of various contour integrals.  In \S\ref{conclude} we discuss the possible implications this has to Real Multiplication, and where this result could be improved.

\section{Holomorphic Theta functions for $L$}\label{holo}

We begin by showing that there are no nontrivial\footnote{A trivial theta function is a nonzero multiple of the exponential function.} holomorphic theta functions for a pseudolattice. 

\begin{prop}\label{nohol}
Let $L$ be a pseudolattice.  There are no nontrivial holomorphic theta functions for $L$.
\end{prop}
\begin{proof}
Note that if $\Theta$ is a theta function for $A_l \in Z^1(L_, \mathcal{H}^\ast)$
and $$B_l(v)=A_l(v)\frac{h(v+l)}{h(v)}$$ for some non-vanishing holomorphic function $h$, then
$\Theta(v)h(v)$ is a theta function for $B_l$.  Hence it suffices to show that there are no nonconstant holomorphic theta functions satisfying \eqref{thetacondition} for a representative of each cohomology class in $Z^1(L,\mathcal{H}^\ast)$. \newline

Suppose $\Theta$ is a holomorphic theta function for a line bundle $\mathcal{L}$.
By Theorem 3.18 of \cite{TaylorII}, the isomorphism class of $\mathcal{L}$ in $H^1(L,\mathcal{H}^\ast)$ has a  unique representative 
$$\mu(l)\hat \sigma(\eta)_l(v)$$
where $\mu \in \textrm{Hom}(L,U(1))$ and $\hat \sigma(\eta)_l(v)$ is defined in (11) of \cite{TaylorII} to be
$$\hat \sigma(\eta)_l(v):=e^{s_\eta{\pi i \over \omega_1}[b^2\omega_2+2bv]} $$
for $l=a \omega_1+b \omega_2$.  Here $L=\omega_1 \B{Z}+\omega_2\B{Z}$, $\eta$ is the image of $\mathcal{L}$ under the Chern class map, and $s_\eta \in \B{Z}$ is the integer associated to $\eta$ under the natural isomorphism $\textrm{Alt}^2(L,\B{Z}) \cong \B{Z}$.  For the details see \cite{TaylorII}.\newline

For $v \in \B{R}$ we have $\left|\mu(l)\hat \sigma(\eta)_l(v) \right|=1$.  Hence for all $v \in \B{R}$ we have 
\begin{equation}\label{periodl}
\left|\Theta(v+l) \right|=\left|\mu(l)\hat \sigma(\eta)_l(v) \Theta(v)\right|=\left| \Theta(v)\right|.
\end{equation}

First note that if $\Theta(v)$ has a zero, then it is identically zero, for the above relation implies that it has an accumulation point of zeros.  Therefore we may assume that $\Theta(v)$ is nonvanishing.\newline

Fix $r \in \B{R}$.  Since $\Theta$ is nonvanishing there exists a function $x_r(v)$ holomorphic in $v$ such that
\begin{equation}\label{s}
{\Theta(v+r) \over \Theta(v)}=e^{2 \pi i x_r(v)}.
\end{equation} 
Without loss of generality we assume that $x_0(v)=0$.  Equation \eqref{periodl} implies that $x_r(v) \in \B{R}$ for all $v \in \B{C}$.  Since $x_r(v)$ is a holomorphic function in $v$ this implies that it is constant. \newline

Now fix $v \in \B{C}$, and let $r,s \in \B{R}$.  
Then 
$${\Theta(v+r+s) \over \Theta(v)}={\Theta(v+r+s) \over \Theta(v+r)}{\Theta(v+r) \over \Theta(v)}. $$
Hence there exists $n(v) \in \B{Z}$ such that
$$
\begin{array}{rcl}
x_{r+s}(v)&=&x_s(v+r)+x_r(v)+2 \pi i n(v)\\
&=&x_s(v)+x_r(v) + 2 \pi i n(v).
\end{array}
$$
Since $x_0(v)=0$ we see that $n(v)=0$, and as a function of $r \in \B{R}$, $x_r(v)$ is a homomorphism. Hence for all $r \in \B{R}$, $x_r(v)=\alpha r$ for some $\alpha \in \B{R}$.\newline

Now consider the left hand side of \eqref{s}.  As $r$ varies over $\B{C}$ this is a holomorphic function.  Hence 
for fixed $v$, there exists a function $x_v(w)$ holomorphic in $w$ such that for all $v \in \B{C}$
$$
{\Theta(v+w) \over \Theta(v)}=e^{2 \pi i x_v(w)}.
$$
Again we assume without loss of generality that $x_v(0)=0$, and therefore $x_z(v)=x_v(z)$ for all $v, w \in \B{C}$. On $\B{R}$ we therefore have $x_v(w)=\alpha w$, and hence by holomorphicity this holds on the whole plane.
\newline

Since $\Theta$ is holomorphic we may compute its derivative along any path.  Let $z \in \B{C}$, and let $\gamma_z(t)$ be the path $z+t$.  Then 
$$\Theta'(z)=\lim_{t \rightarrow 0} {\Theta(\gamma_z(t)) -\Theta(\gamma_z(0)) \over t}=\lim_{t \rightarrow 0} {\Theta(z+t) -\Theta(z) \over t}$$
$$=\lim_{t \rightarrow 0}{e^{2 \pi i \alpha t}-1 \over t}\Theta(z)=2 \pi i \alpha \Theta(z). $$
Hence
$${d \over dz} \log(\Theta(z))=2 \pi i \alpha$$
and hence $\Theta(z)=Ae^{2 \pi i \alpha z}$ for some $A \in \B{C}^\ast$.

\end{proof}

This could be viewed as a set back in defining Real Multiplication analogues to functions which form the foundation of Complex Multiplication, such as the Weierstrass $\wp$-function and modular discriminant $\Delta$.  When $X_\tau$ is the Complex Torus corresponding to the lattice $\Lambda_\tau=\B{Z}+\B{Z}\tau$ there are four holomorphic Jacobi theta functions denoted by $\theta_i(z,\tau)$ for $i=1,2,3,4$.  These are related to the $\wp$ and $\Delta$-functions via the following expressions: 
$$
\wp(z;\tau)=-\log(-\theta_1(z;\tau))''+c  \textrm{  for some constant $c$;}$$
$$
\theta_3(0;\tau)^{24}=\frac{\Delta^2(\frac{\tau+1}{2})}{\Delta(\tau+1)} \textrm{  where $\Delta=\eta^{24}$.}
$$
We may have hoped that the existence of nontrivial holomorphic theta functions associated to Quantum Tori would have enabled us to define similar functions for a real irrational parameter $\theta$ in place of the complex modulus $\tau$.\newline

The nonexistence of nontrivial holomorphic theta functions for pseudolattices leads us to consider the existence of theta functions which are meromorphic. For Complex Tori, elliptic functions and meromorphic theta functions can be constructed out of quotients of holomorphic theta functions \cite{Farkas}.  For Quantum Tori this technique fails due to Proposition \ref{nohol}.  In the next section we examine how it is possible to define meromorphic theta functions for $L$.

\section{Meromorphic Theta Functions for $L$}\label{mero}

Let $\mathcal{H}$ denote the ring of holomorphic functions on $\B{C}$, and denote by $\mathcal{K}$ the field of fractions of $\mathcal{H}$.  Then $\mathcal{K}^\ast$
is the multiplicative group of those meromorphic function which are not identically zero.
Consider the group of 1-cocycles $Z^1(L,\mathcal{K}^\ast)$. These can be
viewed as cocycles corresponding to \emph{meromorphic theta functions} for the pseudolattice $L$.  We saw in the previous section that any holomorphic theta function for $L$ is constant.  This motivates the following question:

\begin{queste}[Existence of meromorphic theta functions for ${L}$]
Does there exist $A_l(v) \in Z^1(L,\mathcal{K}^\ast)$, and a nonconstant meromorphic function $F$ on $\B{C}$ such that for any $l \in L$, $v \in \B{C}$ we have
$$ F(v+l)=A_l(v)F(v)?$$
\end{queste}

\begin{defn}\label{doublesine}
Let $\omega=(\omega_1,\omega_2)$ be a $2$-tuple of elements $\omega_1,\omega_2 \in \B{R}_{>0}$.
The {\it double sine function} with parameter $\omega$ is the unique meromorphic function $S_2^{\omega}(z)$ on $\B{C}$ such that:
\begin{eqnarray}
S_2^{\omega}(z,\omega)=2 \sin \left(\frac{\pi z }{ \omega_2}\right)S_2^{\omega}(z+\omega_1,\omega) \label{1}\\
S_2^{\omega}(z,\omega)=2 \sin \left(\frac{\pi z }{ \omega_1}\right)S_2^{\omega}(z+\omega_2,\omega) \label{2}\\
S_2^{\omega}\left(\frac{\omega_1+\omega_2}{2},\omega_2 \right)=1.
\end{eqnarray}
\end{defn}

This existence of such a function can be deduced from the properties of the double gamma function.  The development of the double gamma function by Barnes in \cite{BarnesII} in 1901 was motivated by Lerch's formula
\begin{equation}\label{lerch}
\log \Gamma(x)=\zeta'(0,x)+{1 \over 2}\log(2 \pi),
\end{equation}
where $\zeta(s,x)$ is the Riemann-Hurwitz zeta function.  For $x,\omega_1,\omega_2 \in \B{R}_{>0}$, and $s \in \B{C}$ with $\Re(x)>0$,  the double  Riemann-Hurwitz zeta function is defined to be
$$\zeta_2(s,x,(\omega_1,\omega_2))=\sum_{n_1,n_2 \in \Bbb{N}}{1 \over (n_1 \omega_1+
n_2 \omega_2 +x)^s}.$$
This series converges absolutely for $\Re(s)>1$ and has an analytic continuation to the complex plane.
The relationship between the gamma function and Riemann-Hurwitz zeta function in \eqref{lerch} motivates the double gamma function $\Gamma_2(x,\omega)$ to be defined by the following relation:
$$\log \left(\Gamma_2(x,\omega)\right)=\left.{\partial \over \partial s}\zeta_2(s,x,\omega)\right|_{s=0}+A$$
where $A$ is some normalising constant.  We can now define the double sine function by the following formula:
$$S_2^\omega(z):={\Gamma_2(\omega_1+\omega_2-z,\omega) \over \Gamma_2(z,\omega)}.$$

\begin{prop}\label{meroauto}
Let $L=\B{Z}\omega_1+\B{Z}\omega_2$ be a pseudolattice.  Then $S_2^\omega(z)$ is a meromorphic theta function for $L$.  More generally, suppose that $G$ is a meromorphic function such that there exist meromorphic functions $f(v)$ and $g(v)$ such that for all $v \in \B{C}$
\begin{eqnarray}
G(v+\omega_1)=f(v)G(v) \label{pd1} \\
G(v+\omega_2)=g(v)G(v). \label{pd2}
\end{eqnarray}
Then $G(v)$ is a meromorphic theta function for $L$.  If $l=n\omega_1+m\omega_2 \in L$ then we have
\begin{equation}\label{periodic}
G(v+l)=A_l(v)G(v)
\end{equation}
where
$$A_l(v):=\prod_{r=0}^{n-1} \prod_{s=0}^{m-1}f(z+r\omega_1)g(z+s\omega_2).$$
\end{prop}
\begin{proof}
It suffices to prove the general case.  It is an immediate consequence of the periodicity relations of \eqref{pd1} and \eqref{pd2} to show that \eqref{periodic} is satisfied.  We need to show that $A_l(v) \in Z^1(L,\mathcal{H}^\ast)$. Let $l_1=n_1 \omega_1+m_1\omega_2$ and $l_2 =n_2 \omega_1+m_2 \omega_2 \in L$, and let $l=l_1+l_2=n\omega_1+m\omega_2$.  Then 
\begin{eqnarray*}
A_{l_1+l_2}(v) & = &
\prod_{r=0}^{n-1} \prod_{s=0}^{m-1}
f(v+r\omega_1)g(v+s\omega_2) \\
&=&\prod_{r=n_2}^{n-1} \prod_{s=m_2}^{m-1}
f(v+r\omega_1)g(v+s\omega_2)\\
& & \qquad  \times \  \prod_{r=0}^{n_2-1} \prod_{s=0}^{m_2-1}
f(v+r\omega_1)g(v+s\omega_2)\\
& = & \prod_{r=0}^{n_1-1} \prod_{s=0}^{m_1-1}
f(v+(r+n_2)\omega_1)g(v+(s+m_2)\omega_2)\\\\
& & \qquad \times \ A_{l_2}(v) \\
& = & A_{l_1}(v+l_2)A_{l_2}(v)
\end{eqnarray*}
\end{proof}

Having exhibited the existence of meromorphic theta functions for pseudolattices, we shall now concern ourselves with their possible application to Real Multiplication.  In the next section we examine the work of Stark and Shintani to Hilbert's twelfth problem, and observe that meromorphic theta functions for pseudolattices have an important role to play in this area.

\section{Stark's Conjecture and Hilbert's Twelfth Problem}   \label{starks conj}
In this section we give a brief overview of a series of conjectures made by Stark concerning the values of $L$-functions associated to number fields at $s=0$.  This leads on to give an account of the work of Shintani, who proved a version of Stark's conjecture in special cases when the ground field was a real quadratic field.  We aim to stress the importance of the double sine function in Shintani's method, and its application in his approach to a solution of Hilbert's twelfth problem for certain real quadratic fields.  Motivated by so called ``higher order'' Stark conjectures, and Shintani's results we will study generalisations of the double sine function in \S\ref{gam}.\newline

\subsection{L-functions and Stark's conjecture}

Let $K$ be a number field, and suppose that $M$ is an abelian extension of $K$ with Galois group $G$.  Class field theory supplies a homomorphism 
$$\tilde \psi_{M/K}: I_K \rightarrow G$$
where $I_K$ denotes the group of fractional ideals of $K$. Let $V$ be a representation of $G$ with character $\chi$. Then define
$$L(\chi,s)=\prod_{\mathfrak{p}} L_\mathfrak{p}(\chi,s)$$
where $\mathfrak{p}$ runs over the prime ideals in $\mathcal{O}_K$ and 
$$L_\mathfrak{p}(\chi,s)=\left(
1-{\chi(\tilde \psi_{M/K}(\mathfrak{p})) N_{K/\Bbb{Q}} \mathfrak{p}^{-s}}\right)^{-1}.$$
Let $S$ be a finite set of places of $K$ which is non-empty and contains all the infinite places of $K$.  We define the L-function associated to $S$ by 
$$L_S(\chi,s):=\prod_{\mathfrak{p} \notin S}L_\mathfrak{p}(\chi,s).$$
These functions are known as L-functions, and when $\chi$ is a nonprincipal character they have analytic continuations to the entire complex plane.  There exists a functional equation for these functions relating their values at $s$ to their values at $1-s$.  \\

We can write
$$L_S(\chi,s)=\sum_{g \in G}\chi(g) L_S(s,g)$$
where 
$$L_S(s,g)=\sum_{
\begin{subarray}{c}
\mathfrak{a}:(\mathfrak{a},S)=1\\
\tilde \psi_{M/K}(\mathfrak{a})=g
\end{subarray}} \frac{1}{N_{K/\B{Q}}(\mathfrak{a})^s}.$$

When $\chi$ is nonprincipal the functional equation implies (see \cite{Tate}) that the order of vanishing of $L_S(\chi,s)$ is equal to 
\begin{equation}\label{rchi}
r(\chi)=\left|\{v \in S: \textrm{$v$ splits completely in $L$} \}\right|.
\end{equation}
Suppose $L$ is ramified at precisely one of the infinite primes, and that $S$ contains precisely the ramified finite primes and the infinite ones.  Then $r(\chi)=1$.\\

In a series of four papers \cite{StarkI,StarkII,StarkIII,StarkIV} between 1971 and 1980, Stark studied the values of the $L$-functions attached to such Galois extensions of number fields at the value $s=1$, which are related via the functional equation to the values at $s=0$.  If as above, the $L$-function has a first order zero at $s=0$, the simple pole of the gamma factor of the functional equation picks out the derivative of $L_S(\chi,s)$ at $s=0$.  Under these conditions studying the value of the $L$-function at $s=1$ is equivalent to studying the value of the derivative at $s=0$.

\begin{conj}[Rank One Abelian Stark Conjecture \cite{Roblot},\cite{Tangedal}]   \label{Stark}
Let $M/K$ be an abelian extension, and $S$ a finite set of places of $K$ containing the infinite ones, one of which splits completely in $M$.  Let $m$ be the number of roots of unity contained in $K$.  There exists an $S$-unit (not necessarily unique) $\varepsilon \in M$
such that for every character $\chi$ of $G$ we have
$$\left.{d \over ds}L_S(\chi,s)\right|_{s=0}=-{1 \over m} \sum_{\sigma \in G}
\chi(\sigma) \log \left|\varepsilon^\sigma \right|_w.$$
\end{conj}

Variations on this conjecture exist for when both infinite primes ramify (known as the Brumer-Stark conjecture), and when $K$ is totally real.  This last case was studied by Tangedal in \cite{Tangedal}.\newline

In the last of Stark's papers he proves a version of Conjecture \ref{Stark} for the cases
case $k=\Bbb{Q}$, and when $k$ is an imaginary quadratic field.  The latter result uses the
the work of Ramachandra in \cite{Ramachandra}, which also was a driving force behind
the work of Shintani, whose work we study in the next section.  \newline

\subsection{Real Quadratic Fields and the work of Shintani}\label{rqfshin}

In 1976 Shintani \cite{ShintaniI} introduced a generalisation of the Riemann-Hurwitz zeta function and proved its analytic continuation to the complex plane.  Shintani used this function to reprove the result of Siegel and Klingen \cite{Hida,ShintaniIII}:
\begin{quote}
Suppose $k$ is a totally real field, and let $\chi$ a character of the ray class group of $F$ modulo an integral ideal $\mathfrak{f}$.   Let $S$ be the finite set of those primes dividing $\mathfrak{f}$.  Then for each $n \in \B{N}$ we have $L_S(1-n, \chi) \in \B{Q}$.
\end{quote}
Shintani showed that for a totally real field $k$, it is possible to express L-functions associated to $k$ as linear combinations of these ``Shintani L-functions'', reducing the study of the value of $L_S(s,\chi)$ to that of  Shintani's L-functions. In a subsequent paper \cite{ShintaniII}, Shintani proved a formula relating the value of his L-functions at $s=1$ to the double gamma function studied by Barnes in \cite{BarnesII}, analogous to the Kronecker limit formula for imaginary quadratic fields.  Using these ideas he went on to prove a refined Stark conjecture for real quadratic fields in \cite{ShintaniIII}. \newline

\subsubsection{Shintani's Limit Formula}
In this section we give an account of Shintani's Kronecker limit formula for real quadratic fields. \newline

Let $F$ be a real quadratic field such that $\textrm{Gal}(F/\Bbb{Q})$ is generated by $\sigma$.
Given an integral ideal $\mathfrak{g}$ of $\mathcal{O}_F$ we let $F_{1,\mathfrak{g}}^+$ denote the group of principal fractional ideals of $F$ generated by those elements $\alpha$ such that
\begin{enumerate}
\item $\alpha$ is totally positive. i.e. $\alpha>0$ and $\alpha^\sigma >0$;
\item $\textrm{ord}_\mathfrak{p}(\alpha-1)>0$ for all $\mathfrak{p} \vert \mathfrak{g}$.
\end{enumerate}
The group $I_F^{\mathfrak{g}}/F_{1,\mathfrak{g}}^+$ is denoted by $G_\mathfrak{g}^+(F)$, and is called the \emph{narrow class group of $F$ modulo $\mathfrak{g}$}, where $I_F^{\mathfrak{g}}$ denotes the group of fractional ideals coprime to $\mathfrak{g}$. When $\mathfrak{g}=\mathcal{O}_F$ we denote this group by $G^+(F)$, and its order by $h^+$.   Given a fractional ideal $\mathfrak{a}$ we let $[\mathfrak{a}]^+$ denote the class it represents in $G^+(F)$.\\

Now fix an integral ideal $\mathfrak{f}$ of $F$, and put
$$S(\mathfrak{f}):=\{\mathfrak{p}: \textrm{$\mathfrak{p}$ is a prime ideal of $\mathcal{O}_F$ dividing $\mathfrak{f}$}\} \cup \{1, \sigma\}.$$
Let $\chi$ be a character of $G_\mathfrak{f}^+(F)$, and suppose $\varepsilon$ is a fundamental totally positive unit of $F$.  \newline

Define a simplicial cone in $\Bbb{R}^2$ by 
$$C:=\{x(1,1)+y(\varepsilon,\varepsilon^\sigma) : x>0, y>0 \}.$$
We choose and fix a set of representatives $\{\mathfrak{a}_1,\mathfrak{a}_2,\ldots ,\mathfrak{a}_{h^+}\}$ of the narrow class group $G^+(F)$ of $F$.  For each $g \in G_\mathfrak{f}^+(F)$ there exists a unique $i$ such that $g=[\mathfrak{fa}_i]^+$ in $G^+(F)$.\newline

With this notation, for $g \in G_\mathfrak{f}^+(F)$ we define the finite set
$$R(g)=\{z=x(1,1)+y(\varepsilon,\varepsilon^\sigma) \in C \cap (\mathfrak{fa}_i)^{-1}:
x \mathfrak{fa}_i \in g, 0<x \leq 1, 0 \leq y <1 \}.$$
Shintani showed that
\begin{equation} \label{partial}
L_S(s,g)=N (\mathfrak{fa}_i)^{-s}\sum_{z=x_1+\varepsilon x_2 \in R(g)}
\zeta(s,(\varepsilon,\varepsilon^\sigma),(x_1,x_2))
\end{equation}
where $\zeta(s,(\varepsilon,\varepsilon^\sigma),(x_1,x_2))$ is a special case of a family of zeta function we will call ``Shintani L-functions''.  Higher dimensional versions of this function were studied in \cite{ShintaniI}, which he used to evaluate the zeta functions associated to totally real algebraic number fields at negative integers, obtaining the result of Siegel and Klingen stated previously.\newline

In \cite{ShintaniII}, Shintani is able to give an explicit formula for the value of the derivative his $L$-function at $s=0$ in terms of Barnes' double gamma function.    Shintani's result can be expressed ass
\begin{equation}  \label{shin}
\left.{d\over ds}
\zeta(s,(\varepsilon,\varepsilon^\sigma),(x_1,x_2))\right|_{s=0}=\log \left(T(x_1+\varepsilon x_2,(\varepsilon,\varepsilon^\sigma))\right) \end{equation}
where
$$T(x_1+\varepsilon x_2,(\varepsilon,\varepsilon^\sigma))=\left\{ {\Gamma _{2}(x_1+x_2\varepsilon
,(1,\varepsilon ))\Gamma _{2}(x+y\varepsilon ^{\sigma },(1,\varepsilon ^{\sigma
}))\over \rho ((1,\varepsilon ))\rho ((1,\varepsilon ^{\sigma }))}
\right\}  e^{{\varepsilon -\varepsilon ^{\sigma }\over 4}
\log \left({\varepsilon ^{\sigma }\over \varepsilon }
\right) \left(x_1^{2}+x_1-{1\over 6}
\right)}.
$$
The numbers $\rho((a_1,a_2))$ are normalising constants which occur in the
theory of the double gamma function \cite{BarnesII}.  The main result of \cite{ShintaniII} is deduced from \eqref{partial} and \eqref{shin}:
\begin{thm}[Shintani, \cite{ShintaniII}]\label{shintaniL}
Let $F$ be a real quadratic field, and $\mathfrak{f}$ an integral ideal of $\mathcal{O}_F$. Let $S=S(\mathfrak{f})$ and suppose $g \in G_\mathfrak{f}^+(F)$.  Then
$$\left.{d\over ds}
L_S(s,g) \right|_{s=0}= \log T(g)$$
where
$$T(g)=\prod_{z=x_1+\varepsilon x_2 \in R(g)} T((z,(\varepsilon, \varepsilon^\sigma)).$$
Hence if $\chi$ is a character of $G_\mathfrak{f}^+(F)$
$$\left.{d \over ds} L_S(\chi,0)\right|_{s=0}=\sum_{g \in G_\mathfrak{f}^+(F)}\chi(g) \log T(g).$$
\end{thm}

This final expression is reminiscent of the one in the Rank One Abelian Stark conjecture (Conjecture \ref{Stark}).  With this comparison, Stark's conjecture suggests that the class invariants $T(g)$ are units in some ray class field over $F$.

\subsection{Shintani's Class Invariants}
In 1978 Shintani produced a paper proving a modified version of Stark's conjecture
for real quadratic fields, subject to various conditions. Astonishingly, he was
unaware of Stark's conjecture when he formulated his results.  \newline

As before, let $F$ be a real quadratic field, and $\mathfrak{f}$ an integral ideal of $F$.
Fix a totally positive integer $\nu$ such that $\nu + 1 \in \mathfrak{f}$, and let $[\nu]_\mathfrak{f}^+$ denote the class it represents in $G_\mathfrak{f}^+(F)$.  By the Existence Theorem of class field theory (Theorem \ref{existe}), there exists an abelian extension $M_\mathfrak{f}$ of $F$ such that the reciprocity map induces an isomorphism
$$G_\mathfrak{f}^+(F) \cong \textrm{Gal}(M_\mathfrak{f}/F) .$$
We shall abuse the notation and shall identify $[\nu]_\mathfrak{f}^+$ with its image under the reciprocity map as an element of this Galois group.  For $g \in \textrm{Gal}(M_\mathfrak{f}/F)$, Shintani studies the value of $L_S(s,g) -L(s,[\nu]_\mathfrak{f}^+ g)$ using
Theorem \ref{shintaniL}.  The properties of $\nu$ imply this has a
particularly nice form:
\begin{equation}  \label{thm1}
L_S(s,g) -L_S(s,[\nu]_\mathfrak{f}^+g)= \sum_{z \in R(g)}
\log \left\{F(z,(1,\varepsilon))F(z^\sigma,(1,\varepsilon^\sigma))
\right\}\end{equation}
where the function $F(z,(1,\varepsilon))$ is related to the double sine function introduced in Definition \ref{doublesine} by
$$  \label{sine}
F(z,(1,\varepsilon)) = S_2^{(1,\varepsilon)}(z)^{-1}.$$
Based on this result Shintani defines the natural class invariant
$$X_\mathfrak{f}(g)=\prod_{z \in R(g)}F(z,(1,\varepsilon))
F(z^\sigma,(1,\varepsilon^\sigma)).$$
With the notation of Theorem \ref{shintaniL} we have 
$$X_\mathfrak{f}(g)=T(g)T([\nu]_\mathfrak{f}^+g)^{-1}.$$
Theorem \ref{shintaniL} implies that if Stark's conjecture is true, the invariants $X_\mathfrak{f}(g)$ should be units.\\

For a subgroup $G$ of $G_\mathfrak{f}^+(F)$, given $c \in G_\mathfrak{f}^+(F)/G$ define
$$X_\mathfrak{f}(c,G)=\prod_{g \in G}X_\mathfrak{f}(cg),$$
and let $M_\mathfrak{f}(G)$ denote the subfield of $M_\mathfrak{f}$ fixed by the elements of $G$.  Using these invariants Shintani proves the following subject to some conditions on $G$ and further rather restrictive hypotheses on the ideal $\mathfrak{f}$.

\begin{thm}   \label{sh}
There exists a positive rational number $m$ such that
\begin{enumerate}
\item The invariant $X_\mathfrak{f}(c,G)^m$ is a unit in the field $M_\mathfrak{f}(G)$.  Moreover for every $g \in G_\mathfrak{f}^+(F)$ we have
$$\left\{X_\mathfrak{f}(c,G)^m \right\}^{\psi_{M_\mathfrak{f}/F}(g)}=X_\mathfrak{f}(cg,G)^m$$
\item Consider the system of invariants
$$\bigcup_{\mathfrak{f}' \parallel \mathfrak{f}}
\left\{
X_{\mathfrak{f}_0}(c,\tilde{G})^m: c \in G_{\mathfrak{f}_0}^+(F)/\tilde{G} \right\}.$$
The union is taken over all divisors $\mathfrak{f}_0$ of $\mathfrak{f}$
which satisfy the same conditions that $\mathfrak{f}$ does, and $\tilde{G}$
is the image of $G$ under the natural homomorphism
$$G_\mathfrak{f}^+(F) \longrightarrow G_{\mathfrak{f}_0}^+(F).$$
Then this system generates $M_\mathfrak{f}(G)$ over $F$.
\end{enumerate}
\end{thm}

The conditions on $G$ imply that precisely one of the infinite primes of $F$
splits in $M_\mathfrak{f}(G)$, so we are in the case considered by the Rank One Abelian Stark conjecture.  Theorem \ref{sh} not only serves to give a special case of Stark's conjecture, but also has the two ingredients listed in the introduction which are required as a solution to Hilbert's twelfth problem: A system of generators with an explicit action of the Galois group.

\section{Generalisations of the Double Gamma Function}\label{gam}

Shintani's results imply that the double sine function will play an important role in any solution to the Rank One Abelian Stark conjecture for real quadratic fields.  The description of this function as a meromorphic theta function for a pseudolattice, and hence a section of a line bundle over a Quantum Torus leads us to question the existence of other such functions.  In this section we generalise the notion of the double gamma function originally defined by Barnes, with a view to investigating its relationship to the values of L-series attached to real quadratic fields.\\

Let $\omega_1,\omega_2 \in \B{R}$ be such that the quotient $\omega_2/\omega_1$
is not negative.  In \cite{BarnesII}, the double gamma function was defined by the integral equation
\begin{equation}\label{gamma}
\Gamma_2(z,\omega):=\exp \left\{{1 \over 2 \pi i} \oint_{I(\lambda,\infty)} e^{-zt}{1 \over (1-e^{-\omega_1t})(1-e^{-\omega_2t})}{\log(-t)+\gamma \over t}dt\right\}.\end{equation}
In this representation and in what follows, for $r \in \B{R}_{>0} \cup \{\infty\}$,  $I(\lambda,r)$ is the contour from $r$ towards zero along the positive real axis to $\lambda$, around zero anticlockwise by a circle of radius $\lambda$ and then out along the real axis to $r$. \newline

We aim to generalise this integral definition to define a family  $\Gamma_2^{r}(z,\omega)$ of functions which satisfy periodicity conditions with respect to the group $\B{Z}\omega_1+\B{Z}\omega_2$, for which we have 
$\Gamma_2^{1}(z,\omega)=\Gamma_2(z,\omega)$.\\

Barnes supplies the following defining relation for the double sine function
$$\log \left(\Gamma_2(z,\omega)\right)=\left.{d \over ds}\zeta_2(s,z,\omega)\right|_{s=0}$$
where
$$\zeta_2(s,z,\omega)=\sum_{n,m=0}^\infty {1 \over (z+m\omega_1+n\omega_2)^s}$$
for $\Re(s)>1$ and $\Re(z)>0$.  The function $\zeta_2(s,z,\omega)$ has meromorphic continuation to the whole plane as a function of $s$ and $z$.

\begin{defn}
Let $r \in \B{N}$, and suppose ${\omega}=(\omega_1,\omega_2) \in \B{R}^2$ is such that the quotient $\omega_2/\omega_1$ is not negative.  For $z \in \B{C}$ define
\begin{equation}\label{generalgam}
\log(\Gamma_2^{r}(z,\omega)):=\left. \left(d \over ds\right)^r \zeta_2(s,z,\omega)\right|_{s=0}.
\end{equation}
\end{defn}

We have an integral formula for $\zeta_2(s,z,\omega)$ given by 
$$\zeta_2(s,z,\omega)={\Gamma(1-s) \over 2 \pi i} \oint_{I(\lambda,\infty)} e^{-zt} {(-t)^{s-1} \over (1-e^{\omega_1t})(1-e^{\omega_2t})}dt.$$
Integrating this $r$ times we obtain an integral expression for $\Gamma_2^r(z,\omega)$:
\begin{equation}\label{matrix}
\log\left( \Gamma_2^r(z,\omega)\right)={1 \over 2 \pi i}\sum_{m=0}^r (-1)^m\binom{r}{m}\Gamma^{(m)}(1) \oint_{I(\lambda,\infty)} {e^{-zt} \over (1-e^{\omega_1t})(1-e^{\omega_2t})}{\log(-t)^{r-m} \over t}dt.
\end{equation}

\begin{defn}
For $\omega=(\omega_1,\omega_2) \in \B{R}^2$ such that the quotient $\omega_2/\omega_1$ is not negative, $z \in \B{C}$ with $\Re(z)>0$ define 
$$
G_2^r(z,\omega)=\exp \left({1 \over 2 \pi i}\oint_{I(\lambda,\infty)} {e^{-zt} \over (1-e^{\omega_1t})(1-e^{\omega_2t})}{\log(-t)^{r} \over t}dt\right)
$$
\end{defn}

\begin{lemma}\label{linG}
Fix $z \in \B{C}$ and $\omega \in \B{R}_{>0}^2$.
Let $W$ denote the field generated over $\B{Q}$ by the values $\Gamma^{(i)}(1)$ for $i=0, \ldots, r$.  Let $V$ vector space over $W$ generated by the values $\log \left(\Gamma_2^j(z,\omega)\right)$ for $j=0,\ldots, r$.  Then $V$ is equal to the vector space over $W$ generated by $\log \left(G_2^j(z,\omega)\right)$ for $i,j=0,\ldots, r$.
\end{lemma}
\begin{proof}
Define a matrix $A$ with coefficients
$$A_{ij}=(-1)^j\binom{i}{j}\Gamma^{(j)}(1) \in W.$$
By \eqref{matrix} we have
$$\log\left(\Gamma^r_2(z,\omega)\right)=\sum_{j=0}^r A_{rj}\log \left(G_2^{r-j}(z,\omega)\right).$$
The matrix $A_{ij}$ is upper triangular, with nonzero diagonal entries, and therefore invertible.
\end{proof}
\begin{cor}
 $G_2^r(z,\omega)$ is a meromorphic theta function on $\B{C}$ for the pseudolattice $L=\B{Z}\omega_1+\B{Z}\omega_2$.
\end{cor}
\begin{proof}
The meromorphicity follows from the meromorphicity of the $\zeta_2(z,s,\omega)$, \eqref{generalgam} and Lemma \ref{linG}.   Observe that
$$\zeta_2(s,z+\omega_1,\omega)={\Gamma(1-s) \over 2 \pi i} \oint_{I(\lambda,\infty)} e^{-zt}[1+(1-e^{\omega_1t})] {(-t)^{s-1} \over (1-e^{\omega_1t})(1-e^{\omega_2t})}dt.
$$
$$=\zeta_2(s,z,\omega)+{\Gamma(1-s) \over 2 \pi i} \oint_{I(\lambda,\infty)} e^{-zt} {(-t)^{s-1} \over (1-e^{\omega_2t})}dt.
$$
The second term is equal to a zeta function $\zeta_1(s,z,\omega_2)$ which has meromorphic continuation to the whole plane \cite{BarnesI}.  Let $\Gamma_1^r(z,\omega_2)=\exp\left(\zeta_1^{(r)}(0,z,\omega_2) \right)$, and hence 
$$\Gamma_2^r(z+\omega_1,\omega)=\Gamma_1^r(z,\omega_1)\Gamma_2^r(z,\omega).$$
A similar expression holds for $\Gamma_2(z+\omega_2,\omega)$.\\

Hence the functions $\Gamma_2^r(z,\omega)$ are meromorphic theta functions for $L$.  By Lemma \ref{linG}, the functions $G_2^r(z,\omega)$ are.
\end{proof}

Our final aim is to write the higher derivatives of L-functions associated to real quadratic fields in terms of meromorphic theta functions for a pseudolattice.  In order to achieve this we will need to introduce another function, which does not seem to have any analogy in Shintani's work.\newline

For $t,u,z,v \in \B{C}$ and ${\omega}, {\lambda} \in \B{R}^2$ define
\begin{equation}\label{g}
\mathfrak{g}(t,u,z,v,{\omega},{\lambda})= {e^{zt}e^{(\left| \lambda \right|-v)tu} \over (1-e^{t(\omega_1+u \lambda_1)})(1-e^{t(\omega_2+u \lambda_2)})}-{e^{zt} \over (1-e^{t\omega_1})(1-e^{t \omega_2})}
\end{equation}
where $\left|\lambda \right|=\lambda_1+\lambda_2$.  This is a holomorphic function in $u$ with a zero at $u=0$.  We define a family of functions $C_N(t,v,\omega,\lambda)$ indexed by $N \in \B{N}$ by
\begin{equation}\label{coeffg}
\mathfrak{g}(t,u,z,v,{\omega},{\lambda})+{e^{zt} \over (1-e^{t\omega_1})(1-e^{t \omega_2})}=\sum_{N=0}^\infty e^{tz}C_N(t,v,{\omega},{\lambda})u^N. 
\end{equation}
We note that
$$C_0(t,v,{\omega},{\lambda})={1 \over (1-e^{t\omega_1})(1-e^{t \omega_2})},$$
and hence 
$$\log(G_2^r(z,{\omega}))={1 \over 2 \pi i}\oint_{I(\lambda,\infty)} {\log(-t)^r \over t} e^{(\left| {\omega} \right|-z)t}C_0(t,v,{\omega},{\lambda}) dt$$
for any ${\lambda} \in \B{R}^2$, $v \in \B{C}$ where $\left|{\omega} \right|=\omega_1+\omega_2$.\\

Suppose $h$ is a function in a real variable vanishing at $0$.
Let $J$ be the operator defined on such a function by
$$J(h)(u):=-{1 \over 2 \pi i}\int_0^u {1 \over t} h(t) dt.$$
Let $\mathfrak{g}(u)$ be the function defined in \eqref{g}, considered as a function of $u$.
Then we have
$$J^{k}(\mathfrak{g}(u))(1)=(-1)^k\sum_{N=1}^{\infty} { C_{N}(t,v,\omega,\lambda) \over N^{k}},$$
where the functions $C_{N}(t,v,\omega,\lambda)$ are as defined in \eqref{coeffg}.  Note that this can be viewed as a variety of zeta function.
\begin{defn}
Suppose ${\omega}, {\lambda} \in \B{R}^2$ are such that neither of the quotients $\omega_2/\omega_1$ or $\lambda_2/\lambda_1$ are negative.  For $z,v \in \B{C}$ and $q,k \in \B{N}$ define
\begin{equation}\label{generalh}
H^{q,k}(z,v,{\omega},{\lambda}):={1 \over 2 \pi i}\oint_{I(\lambda,\infty)} e^{(\left| \omega \right|-z)t}J^{k}(\mathfrak{g}(u))(1){ \log(-t)^q \over t}dt.
\end{equation}
\end{defn}
\begin{prop}
For all $k,q, \in \B{N}$, the integral of \eqref{generalh} converges for $\Re(z)>S$ for some $S$ depending on $v$, $\omega$ and $\lambda$.  In this region the integral defines an analytic function $H^{k,q}(z,v,\omega,\lambda)$, which is a theta function in $z$ for the pseudolattice $L=\B{Z}\omega_1+\B{Z}\omega_2$.
\end{prop}
\begin{proof}
The proof of this result is the subject of \S\ref{analH}.
\end{proof}
This last result follows as a result of some crude estimates using Cauchy's integral formula for the derivative of a holomorphic function. We conjecture that this may be strengthened:
\begin{conj}\label{conjh}
For all $k,q \in \B{N}$, $v \in \B{C}$ and $\omega,\lambda \in \B{R}_{>0}^2$ the integral in \eqref{generalh} defines a meromorphic function $H^{q,k}(z,v,\omega,\lambda)$, which as a function in $z$ is a theta function for $L=\B{Z}\omega_1+\B{Z}\omega_2$.
\end{conj}

We will implicitly assume Conjecture \ref{conjh} for the remainder of this thesis.

\section{The derivative of L-functions of real quadratic fields}\label{der}

Let $F$ be a real quadratic field and suppose $\mathfrak{f}$ is an integral ideal of $F$. Let $S$ be a finite set of primes of $F$ containing those primes dividing $\mathfrak{f}$.  Let $\chi$ be a character of the group $G_\mathfrak{f}^+(F)$, and let $L_S(\chi,s)$ denote the corresponding L-function.  In \S\ref{rqfshin}, we saw from the results of Shintani that under certain circumstances we can write the value of $L_S'(\chi,0)$ as a linear combination of special values of meromorphic theta functions for pseudolattices lying in $F$.  In this section will prove the following
\begin{thm}\label{special}
Let $F$ be a real quadratic field, $\mathfrak{f}$ an integral ideal of $F$, and $\chi$ a character of $I_F^\mathfrak{f}$.  Let $m \in \B{N}$, and let $L^{(m)}(\chi,s)$ denote the $m^{th}$ derivative of the L-function with respect to $s$.  We may write $L_S^{(m)}(\chi,0)$ as an element of the field $K^m_\mathfrak{f}(F)$ generated over $F$ by 
\begin{enumerate}
\item \label{trans} $2 \pi i$, the values $\Gamma^{(j)}(1)$ for $j=0, \ldots, m$.  The maximal power of $2 \pi i $ which occurs is $m+1$;
\item \label{roots}the roots of unity of order $p$, where $p$ is the maximal order of an element of $G_\mathfrak{f}^+(F)$;
\item \label{logs} the logarithms of a finite number of elements $N_i \in F$ (which are specified in the statement of Lemma \ref{Lfunc});
\item \label{polylogs} the values $Li_n\left(-\varepsilon^\sigma/\varepsilon \right)$, $Li_n\left(-\varepsilon/\varepsilon^\sigma \right)$ and $Li_n\left(-1\right)$ for $n=1 \ldots m+1$, where $Li_n$ denotes the $n^{th}$ polylogarithm function, and $\varepsilon$ is a generator for the group of totally positive units of $F$;
\item \label{spec} the special values 
$$
\begin{array}{c}
\log(G_2^{r}(x_1^i+\varepsilon x_2^i,(1,\varepsilon))) \\  
\log(G_2^{r}(x_1^i+\varepsilon^\sigma x_2^i,(1,\varepsilon^\sigma)))\\
\log(H^{r,k}(x_1^i+\varepsilon x_2^i,x_1^i+\varepsilon^\sigma x_2^i,(1,\varepsilon),(1,\varepsilon^\sigma)))\\
\log(H^{r,k}(x_1^i+\varepsilon^\sigma x_2^i,x_1^i+\varepsilon x_2^i,(1,\varepsilon^\sigma),(1,\varepsilon)))
\end{array}$$
where $\underline{x}^i$ is one of a finite set of pairs of element of $F$ determined by $F$ and the choice of $\varepsilon$. The highest value of $r$ and $k$ which occurs is $m$.
\end{enumerate}.
\end{thm}
\begin{rmk}
Throughout the proof I will refer to fields generated over $\B{Q}$ or $F$ by some combinations of these generators.  For example, if I wish to refer to the field generated over $\B{Q}$ by those elements in statements \ref{roots}, \ref{polylogs} and \ref{spec} in the statement of Theorem \ref{special}, I shall denote this field by $\B{Q}([\ref{roots}],[\ref{polylogs}], [\ref{spec}]).$
\end{rmk}

We will break the proof up in to several stages.  The first stage is to recall that we can write the L-function of $F$ as a finite sum of ``Shintani L-functions''.
\begin{defn}[Shintani L-function]
Let $\underline{a}=(a_1,a_2),\underline{x}=(x_1,x_2) \in \B{R}^2$.  Then we define the Shintani L-function $\zeta(s,\underline{a},\underline{x})$ for 
$\Re(s)>1$ by 
\begin{equation}\label{shintani}
\zeta(s,\underline{a},\underline{x})=\sum_{m,n=0}^\infty {1 \over (x_1+m+(x_2+n)a_1)^s(x_1+m+(x_2+n)a_2)^s}.
\end{equation}
\end{defn}
Elements of the proof of the following result were discussed in \S\ref{rqfshin} when we discussed Shintani's Limit Formula:
\begin{lemma}[Shintani, \cite{ShintaniII}] \label{Lfunc}
Let $\varepsilon>1$ be a generator for the group of totally positive units of $F$, and let $\sigma$ be the non trivial element of $Gal(F/\B{Q})$.
There exists $N \in \B{N}$, 2-tuples $\underline{x}_1, \ldots, \underline{x}_N \in F^2$, elements
$N_i \in F$ and $c_i \in \mu_p$ such that 
\begin{equation}\label{L}
L_S(\chi,s)=\sum_{i=1}^N c_i N_i^s \zeta(s,(\varepsilon, \varepsilon^\sigma),\underline{x}_i).
\end{equation}
\end{lemma}
Differentiating the expression for $L_S(\chi,s)$ in \eqref{L} $m$ times with respect to $s$, we see that at $s=0$ the derivative of the $L$-function is given by
\begin{equation}\label{Lder}
L_F^{(m)}(0,\chi)=\sum_{j=0}^m \sum_{i=1}^N c_i\binom{m}{j}  \log(N_i)^{m-j} \zeta^{(j)}(0,(\varepsilon, \varepsilon^\sigma),\underline{x}_i) .
\end{equation}
This expression shows the need to adjoin the roots of unity $\mu_p$ and the values $\log(N_i)$ which are mentioned in parts \ref{roots} and \ref{logs} of the statement of Theorem \ref{special}.  With this result in mind, Theorem \ref{special} will follow if we can prove
the following:
\begin{prop} \label{special2}
Let $m \in \B{N}$ and suppose $\underline{x} \in F^2$. Then with the notation of Theorem \ref{special}, $\zeta^{(m)}(0,(\varepsilon,\varepsilon^\sigma), \underline{x}) \in K^m_\mathfrak{f}(F)$.
\end{prop}

\section{Proof of Theorem \ref{special}}\label{proof}
We will prove Theorem \ref{special} by proving Proposition \ref{special2}.\newline

An integral formula for $\zeta(s,\underline{a},\underline{x})$ is given 
in \cite{ShintaniII} as
\begin{equation}\label{zeta}
4 \pi^2 {(1+e^{2 \pi i s}) \over \Gamma(1-s)^2} \zeta(s,a,x)=
\int_{I(\lambda,\infty)} (-t)^{2s}{dt \over t}
\int_{I(\lambda,1)} u^s {du \over u}
[g(t,tu)+g(tu,t)]
\end{equation}
where 
$$g(t_1,t_2)={e^{(1-x_1)(t_1+t_2)+(1-x_2)(a_1t_1+a_2t_2)} \over (1-e^{t_1+t_2})(1-e^{a_1t_1+a_2t_2})}.$$

We shall proceed by induction.  

\begin{prop}
For $m=0,1$, $\zeta^{(m)}(0,(\varepsilon,\varepsilon^\sigma), \underline{x}) \in K^m_\mathfrak{f}(F)$.
\end{prop}
\begin{proof}
These results follow from the statement and proof of Proposition 3 of \cite{ShintaniII}.  We let $B_1$ and $B_2$ denote the first and second Bernoulli polynomials, which have coefficients in $\B{Q}$.  The statement of this result implies that
$$\zeta^{(1)}(0,(\varepsilon,\varepsilon^\sigma),\underline{x})=\log\left( \Gamma_2^1(x_1+x_2\varepsilon,x_1+x_2\varepsilon^\sigma,(1,\varepsilon),(1,\varepsilon^\sigma))\right) $$
$$+\log\left(\Gamma_2^1(x_1+x_2\varepsilon^\sigma,x_1+x_2\varepsilon,(1,\varepsilon^\sigma),(1,\varepsilon))\right) +
{\varepsilon^\sigma-\varepsilon \over 4 \varepsilon \varepsilon^\sigma}
\log\left({\varepsilon^\sigma  \over \varepsilon}\right)B_2(x_1).$$
We may rewrite the final term as
$${\varepsilon^\sigma-\varepsilon \over 4 \varepsilon \varepsilon^\sigma}\log\left({\varepsilon^\sigma  \over \varepsilon}\right)B_2(x_1)={\varepsilon^\sigma-\varepsilon \over 4\varepsilon \varepsilon^\sigma}\left[Li_1(-\varepsilon^\sigma/\varepsilon)-Li_1(-\varepsilon/\varepsilon^\sigma) \right]B_2(x_1)$$
since $Li_1(x)=-\log(1-x)$.\newline

In the proof of this result, Shintani also shows that 
$$\zeta(0,(\varepsilon,\varepsilon^\sigma), \underline{x})={1 \over 4}\left({1 \over \varepsilon}+{1 \over \varepsilon^\sigma} \right)B_2(x_1)+B_1(x_1)B_1(x_2)+{1 \over 4}(\varepsilon+\varepsilon^\sigma)B_2(x_2).$$ 
Hence we may write the null values of these derivatives of the zeta function in terms of the double gamma function.  Since the field $W$ of Lemma \ref{linG} is contained in $\B{Q}([\ref{trans}])$, the result follows.
\end{proof}

Now fix $m \in \B{N}$, and assume the inductive hypothesis holds for all values of $b$ less than $m$:

\begin{quote}
If $\underline{x} \in F^2$ then $\zeta^{(b)}(0,(\varepsilon,\varepsilon^\sigma), \underline{x}) \in K^m_\mathfrak{f}(F)$ for all $b=0 \ldots m-1$. 
\end{quote}
We need to show that $\zeta^{(m)}(0,(\varepsilon,\varepsilon^\sigma), \underline{x}) \in K^m_\mathfrak{f}(F)$. \newline

Differentiate both sides of \eqref{zeta} $m$ times with respect to $s$.  Evaluating at $s=0$ we see that the left hand side of the result is a finite sum of terms of the form
\begin{equation}\label{LHS}
T_{a,b}^m:=2 \pi^2(2 \pi i)^a \zeta^{(b)}(0,\underline{a},\underline{x})\left.\left(d \over ds \right)^{m-a-b} \Gamma(1-s)\right|_{s=0}\end{equation}
for $a,b \in \B{N}$ such that $a+b \leq m$.  By our inductive hypothesis, if $b \neq m$ then $T_{a,b}^m \in K^m_\mathfrak{f}(F)$.  Note that it is at this point we are required to adjoin the higher derivatives of the gamma function in the statement of Theorem \ref{special}.  To prove that  $\zeta^{(m)}(0,(\varepsilon,\varepsilon^\sigma), \underline{x}) \in K^m_\mathfrak{f}(F)$ it is therefore sufficient to show that the $m^{th}$ derivative of the right hand side of \eqref{zeta} lies in $K^m_\mathfrak{f}(F)$. The rest of this section is devoted to proving that this is indeed true.\newline

When we differentiate $m$ times with respect to $s$, the right hand side of \eqref{zeta} becomes
$$I_m(s):=
\int_{I(\lambda,\infty)} (-t)^{2s}{dt \over t}
\int_{I(\lambda,1)} u^s {du \over u}
[g(t,tu)+g(tu,t)][2 \log(-t)+ \log (u) ]^m.$$
Using the binomial theorem we see that at $s=0$
\begin{equation}\label{split}
I_m(0)=\sum_{p=0}^m \binom{m}{p} 2^p I_{p,m-p}
\end{equation}
where
$$I_{p,q}:=\int_{I(\lambda,\infty)}{dt \over t}
\int_{I(\lambda,1)}  {du \over u}
[g(t,tu)+g(tu,t)]\log(-t)^p\log (u)^q.$$
In order to evaluate integrals of this form we will first consider the integrals
$$A_{p,q}:=\int_{I(\lambda,\infty)}{dt \over t}
\int_{I(\lambda,1)}  {du \over u}
g(t,tu)\log(-t)^p\log (u)^q$$
where $p=m-q$ and $q=0,1,\ldots, m$.  Once we have evaluated the integrals $A_{p,q}$, we shall be able to use our result to evaluate the integrals
\begin{equation}\label{B}
B_{p,q}:=\int_{I(\lambda,\infty)}{dt \over t}
\int_{I(\lambda,1)}  {du \over u}
g(tu,t)\log(-t)^p\log (u)^q.
\end{equation}
Having shown that both $A_{p,q}$ and $B_{p,q}$ lie in $K^m_\mathfrak{f}(K)$, since $I_{p,q}=A_{p,q}+B_{p,q}$ it shall follow that $I_{p,q} \in K^m_\mathfrak{f}(F)$. Hence by \eqref{split} we will have shown $I_m(0) \in K^m_\mathfrak{f}(F)$.  \newline

We begin by noting that 
$$
\begin{array}{rcl}
g(t,tu)&=& \mathfrak{g}(t,u,(1-x_1)+a_1(1-x_2),(1-x_1)+a_2(1-x_2),(1,a_1),(1,a_2))\\
&=&-{e^{zt} \over (1-e^{t})(1-e^{a_1t})}+\sum_{N=1}^\infty e^{tz}C_N(t,v,(1,a_1), (1,a_2))u^N.
\end{array}$$
where $z:=(1-x_1)+a_1(1-x_2)$ and $v=(1-x_1)+a_2(1-x_2)$.\newline

Before we proceed we make a remark which will simplify our calculations.  Note that the integral expression we have for $\zeta(s,a,\underline{x})$ in \eqref{zeta} is independent of $\lambda$ for sufficiently small $\lambda$, and hence $I_m$ is independent of $\lambda$.  Hence $I_m=\lim_{\lambda \rightarrow 0} I_m$. Now suppose we can write $I_m$ as the sum of finitely many integrals:
$$I_m=\sum_{i=1}^n \int_{I(\lambda,\infty)} \int_{I(\lambda,1)} f_i(t,u) du \ dt.$$
Then providing each of the limits is finite, we have  
$$I_m=\lim_{\lambda \rightarrow 0} I_m=\sum_{i=1}^n \lim_{\lambda \rightarrow 0}\int_{I(\lambda,\infty)} \int_{I(\lambda,1)} f_i(t,u) du \ dt.$$
We will use this idea to calculate the integrals $A_{p,q}$.  There are three cases to consider: 
\begin{itemize}
\item When $q=0$ we have 
$$A_{p,0}=-\int_{I(\lambda,\infty)}{dt \over t}
\int_{I(\lambda,1)}  {du \over u}
 e^{zt}C_0(t,v,(1,a_1), (1,a_2))\log(-t)^p$$
$$+\int_{I(\lambda,\infty)}{dt \over t}
\int_{I(\lambda,1)}  {du \over u}\sum_{N=1}^\infty e^{tz}C_N(t,v,(1,a_1), (1,a_2))u^N\log(-t)^p.$$
The second integral vanishes since the integrand does not have poles on or within the contour traced as $u$ traces $I(\lambda,1)$.  By our definitions 
$$
\begin{array}{rcl}
A_{p,0}&=&4 \pi^2 \log\left(G_2^p(1+a_1-z,(1,a_1)) \right)\\
&=&4 \pi^2 \log\left(G_2^p(x_1+a_1x_2,(1,a_1)) \right)
\end{array}.
$$

\item Now consider the case when neither $p$ or $q$ are zero.  Then 
$$A_{p,q}=-\int_{I(\lambda,\infty)}{dt \over t}
\int_{I(\lambda,1)}  {du \over u}
 {e^{zt} \over (1-e^{t})(1-e^{a_1t})}\log(-t)^p\log (u)^q$$
 $$+\int_{I(\lambda,\infty)}{dt \over t}
\int_{I(\lambda,1)}  {du \over u}\sum_{N=1}^\infty e^{tz}C_N(t,v,(1,a_1), (1,a_2))u^N\log(-t)^p\log (u)^q.$$
Using Lemma \ref{log1} (found in the Appendix) we find that this is equal to 
$$-{(2 \pi i)^{q+1} \over q+1}\int_{I(\lambda,\infty)}{dt \over t}
{e^{zt} \over (1-e^{t})(1-e^{a_1t})}\log(-t)^p$$
 $$+\sum_{k=1}^{q-1} (-1)^k {q! \over (q-k)!} (2 \pi i)^{q-k} \int_{I(\lambda,\infty)}{dt \over t}
\sum_{N=0}^\infty e^{tz}{C_N(t,v,(1,a_1), (1,a_2)) \over N^{k+1}}\log(-t)^p.$$
Using the definitions of \S\ref{gam} we find that
$$A_{p,q}=-{(2 \pi i)^{q+2} \over q+1} \log \left(G_2^p(x_1+a_1x_2,(1,a_1))\right)$$
$$-(2 \pi i)^{q+1}\sum_{k=0}^{q-1}{q! \over (q-k)!} \log \left(H^{p,k+1}(x_1+a_1x_2,x_1+a_2x_2,(1,a_1),(1,a_2))\right). $$

\item Finally we consider the case when $p=0$.  In this case 
$$A_{0,q}=-\int_{I(\lambda,\infty)}{dt \over t}
\int_{I(\lambda,1)}  {du \over u}
 {e^{zt} \over (1-e^{t})(1-e^{a_1t})}\log (u)^q$$
 $$+\int_{I(\lambda,\infty)}{dt \over t}
\int_{I(\lambda,1)}  {du \over u}\sum_{N=1}^\infty e^{tz}C_N(t,v,(1,a_1), (1,a_2))u^N\log (u)^q.$$
The coefficient of $t^{-1}$ in the integrand of the first integral is calculated to be
$${1 \over 2}{1+u \over a_1+a_2u} B_2(x_1)+B_1(x_1)B_1(x_2)+{1 \over 2}{a_1+a_2u \over 1+u} B_2(x_2).$$
We use Lemma \ref{log1} again to calculate the integral over $I(\lambda,1)$ to find that 
$$A_{0,q}=-\int_{I(\lambda,1)}\left[{1 \over 2}{1+u \over a_1+a_2u} B_2(x_1)+B_1(x_1)B_1(x_2)+{1 \over 2}{a_1+a_2u \over 1+u} B_2(x_2) \right] {\log (u)^q \over u}$$
 $$+\sum_{k=0}^{q-1} (-1)^k {q! \over (q-k)!} (2 \pi i)^{q-k}\int_{I(\lambda,\infty)}{dt \over t}
\sum_{N=1}^\infty e^{tz}{C_N(t,v,(1,a_1), (1,a_2)) \over N^{k+1}}.$$
Using the definitions of \S\ref{gam} this simplifies to 
$$A_{0,q}=-\int_{I(\lambda,1)}\left[{1 \over 2}{1+u \over a_1+a_2u} B_2(x_1)+B_1(x_1)B_1(x_2)+{1 \over 2}{a_1+a_2u \over 1+u} B_2(x_2) \right] {\log (u)^q \over u}$$
$$-(2 \pi i)^{q+1}\sum_{k=0}^{q-1}{q! \over (q-k)!} \log \left(H^{0,k+1}(x_1+a_1x_2,x_1+a_2x_2,(1,a_1),(1,a_2))\right).$$   
To prove the result it suffices to show that the first integral lies in $K^m_\mathfrak{f}(F)$.  Note that
$$\int_{I(\lambda,1)}{\log(u)^q \over u}{1+u \over a_1+a_2u} du \qquad \qquad \qquad \qquad \qquad \qquad \qquad \qquad \qquad $$
$$\qquad ={1 \over a_1a_2} \int_{I(\lambda,1)} \log(u)^q \left[ {1 \over u}-{a_2 \over a_1+a_2u} \right]\left[(1-a_1)+(a_1+a_2u)\right]du. $$
We are therefore reduced to calculating the following integrals
\begin{equation}\label{m0}
\int_{I(\lambda,1)} \log(u)^q du;
\end{equation}
\begin{equation}\label{m=-1}
\int_{I(\lambda,1)} {\log(u)^q \over u} du;
\end{equation}
\begin{equation}\label{poly}
\int_{I(\lambda,1)} {\log(u)^q \over a_1+a_2u}du.
\end{equation}
Lemma \ref{log1} shows that both the integrals \eqref{m0} and \eqref{m=-1} lie in $K^m_\mathfrak{f}(F)$.  By the remark made earlier, to determine \eqref{poly} it is sufficient to evaluate
$$\label{m=0}
\lim_{\lambda \rightarrow 0}\int_{I(\lambda,1)} {\log(u)^q \over a_1+a_2u}du.
$$
It is easy to see that the integral around the circular path is $O(\lambda)$, so tends to $0$ as $\lambda \rightarrow 0$.  This reduces the evaluation of \eqref{poly} to that of 
$$
\lim_{\lambda \rightarrow 0} \left\{\int_\lambda^1 {(\log(u)+2 \pi i)^q \over a_1+a_2u} du- \int_\lambda^1 {\log(u)^q \over a_1+a_2u} du\right\}.
$$ 
Expanding this using the binomial theorem we are reduced to showing that the following expression lies in $K^m_\mathfrak{f}(F)$:
$$
\lim_{\lambda \rightarrow 0} \int_\lambda^1 {\log(u)^q \over a_1+a_2u} du.
$$
This is proved in Lemma \ref{polylog}.
\end{itemize} 
Hence we have shown that for all $p=0 \ldots m$, $A_{p,m-p} \in K_\mathfrak{f}(F)$.  The details of the proof describes this in more detail:

\begin{lemma}\label{a}
There exist $\alpha_l, \beta_{r,k} \in F([\ref{trans}])$ and polynomials $Q_1,Q_2$ and $Q_3 \in \B{Q}([\ref{trans}],[\ref{polylogs}])[X_1,X_2]$ such that 
$$A_{p,q}=\sum_l \alpha_l \log\left(G_2^l(x_1+a_1x_2,(1,a_1))\right)\qquad \qquad \qquad \qquad $$
$$\qquad +\sum_{r,k} \beta_{r,k} \log\left(H^{r,k}(x_1+a_1x_2,x_1+a_2x_2,(1,a_1),(1,a_2))\right)$$
$$+B_2(x_1)Q_1[a_1,a_2]+B_1(x_1)B_1(x_2)Q_2[a_1,a_2]+B_2(x_2)Q_3[a_1,a_2].$$
\end{lemma}
In order to show that $I_m(0) \in K^m_\mathfrak{f}(F)$ we need to show that $B_{p,m-p} \in K^m_\mathfrak{f}(F)$ where $B_{p,q}$ is defined in \eqref{B}.  However, note that 
$$
\begin{array}{rcl}
g(tu,t)&=& \mathfrak{g}(t,u,(1-x_1)+a_2(1-x_2),(1-x_1)+a_1(1-x_2),(1,a_2),(1,a_1))\\
&=& \mathfrak{g}(t,u,v,z,(1,a_2),(1,a_1))\\
\end{array}
$$
If $A_{p,q}$ has an explicit description as in the statement of Lemma \ref{a}, then the above expression implies that we have 
$$B_{p,q}=\sum_l \alpha_l \log\left(G_2^l(x_1+a_2x_2,(1,a_2))\right)\qquad \qquad \qquad \qquad $$
$$\qquad +\sum_{r,k} \beta_{r,k} \log\left(H^{r,k}(x_1+a_2x_2,x_1+a_1x_2,(1,a_2),(1,a_1))\right)$$
$$+B_2(x_1)Q_1[a_2,a_1]+B_1(x_1)B_1(x_2)Q_2[a_2,a_1]+B_2(x_2)Q_3[a_2,a_1].$$

Hence $B_{p,q} \in K^m_\mathfrak{f}(F)$, and therefore $I_{p,q}=A_{p,q}+B_{p,q}$ does.  By \eqref{split} $I_m(0) \in K^m_\mathfrak{f}(F)$. 
\begin{flushright}$\Box$ \end{flushright}

\section{Meromorphic Theta functions and Stark's conjecture}\label{conclude}

Let $F$ be a real quadratic number field, and let $\mathfrak{f}$ be an integral ideal of $F$.  Let $T$ be a finite set of primes, such that the finite primes in $T$ are precisely those dividing $\mathfrak{f}$, and suppose $\chi$ is a character defined on $I_F^\mathfrak{f}$.  Let $S$ be a set of places of $F$ containing $S$.  Then the relationship between the L-functions corresponding to $T$ and $S$ is given by 
\begin{equation}\label{T}
L_S(\chi,s)=L_T(\chi,s)\prod_{\mathfrak{p} \in T \setminus S}L_{\mathfrak{p}}(\chi,s) .
\end{equation} 
In \S\ref{starks conj} we examined the Rank One Abelian Stark conjecture, which concerned the case when the L-function had a simple zero at $s=0$.  Higher order conjectures, exist when the L-function has zeros of order greater than one.  These conjectures link the values of derivatives of L-functions over number fields to the existence of units in Galois extensions, and have been the subject of study by Tate \cite{Tate} and Rubin \cite{Rubin}.  In the light of these conjectures we suggest that information concerning the derivatives of L-functions associated to real quadratic fields  will be relevant to Real Multiplication.\\

In \S\ref{der} we identified a certain field $K_\mathfrak{f}^m(F)$ in which the value $L_S^{(m)}(\chi,0)$ lies.  If $S$ contained both real primes of $F$ then we observe from \eqref{rchi} that $r(\chi)=2$, and $L_S(\chi,0)$ has a second order zero at $s=0$.  Although we have not explicitly done so, using our method it is possible to give an exact formula for the value $L^{(2)}(\chi,0)$ in terms of the theta functions $G_2^r$ and $H^{k,r}$, and hence to formulate a conjecture on how these functions may define units in a class field above $F$.\\

When $T$ is a set of primes containing $S$, the L-function $L_T(\chi,0)$ may have zeros at $s=0$ of arbitrary order.  Using \eqref{T} we see that 
$$L_T^{(m)}(\chi,0) \in K_\mathfrak{f}^m(F)\left(\{\log(N_{F/\B{Q}}(\mathfrak{p})): \mathfrak{p} \in T \setminus S\} \right).$$
Using analogous techniques to the proof of Theorem \ref{special} we could obtain an explicit formula for this value in terms of meromorphic theta functions for pseudolattices in $F$.\\

Theorem \ref{special} is a far cry from an immediate application to Hilbert's twelfth problem.  It does not give an explicit description for the value of the derivative of the L-function (although this is implicit in the proof), and the field $K_\mathfrak{f}^m(F)$ is clearly not a number field.  However, motivated by a technique of Shintani's \cite{ShintaniIII} we can write certain L-values purely in terms of theta functions, without the transcendental constant terms:

\begin{thm}\label{specialone}
Let $F$ be a real quadratic field, $\mathfrak{f}$ an integral ideal of $F$, and $\chi$ a character of $G_\mathfrak{f}(F)$.  Let $m,n \in \B{N}$.  Let $\mu \in F$ be a totally positive element of $\mathcal{O}_F$ such that 
$\nu \equiv 1 \mod{\mathfrak{f}}$, and suppose $g \in G_\mathfrak{f}^+(F)$.  If $S=S(\mathfrak{f})$ then the value $L_S^{(m)}(0,g)-L_S^{(m)}(0,[\nu]_\mathfrak{f}^+g)$ is an element of the field $K_\mathfrak{f}(F)$ generated over $F$ by 
a finite number of meromorphic theta functions.
\end{thm}
\begin{proof}
By the proof of Theorem 1 of \cite{ShintaniIII}, with the notation used previously we have
$$L_S(s,g)=N(\mathfrak{fa}_j)^{-s} \sum_{x_1+\varepsilon x_2 \in R(g)} \zeta(s,(\varepsilon,\varepsilon^\sigma),(x_1,x_2)).$$
Differentiating $m$ times with respect to $s$ we get
$$L_S^{(m)}(0,g)=\sum_{x_1+\varepsilon x_2 \in R(g)}\sum_{k=0}^m \binom{m}{k}\log(N(\mathfrak{fa}_j))^{m-k} \zeta^{(k)}(0,(\varepsilon,\varepsilon^\sigma),(x_1,x_2)).$$
By Lemma \ref{a}, by employing an induction argument it is easy to show that there exist coefficients $a_l,b_{r,k} \in F([\ref{trans}])$ and polynomials $R_1,R_2$ and $R_3 \in \B{Q}([\ref{trans},\ref{polylogs}])$ such that
$$\zeta^{(m)}(0,(\varepsilon,\varepsilon^\sigma),(x_1,x_2))=\sum_l a_l \log\left(G_2^l(x_1+x_2\varepsilon,(1,\varepsilon))\right)\qquad \qquad \qquad \qquad $$
$$\qquad +\sum_l a_l \log \left(G_2^l(x_1+x_2\varepsilon^\sigma,(1,\varepsilon^\sigma)) \right)\qquad \qquad \qquad \qquad $$
$$\qquad +\sum_{r,k} b_{r,k} \log\left(H^{r,k}(x_1+x_2\varepsilon,x_1+x_2\varepsilon^\sigma,(1,\varepsilon),(1,\varepsilon^\sigma))\right)$$
$$\qquad +\sum_{r,k} b_{r,k} \log\left(H^{r,k}(x_1+x_2\varepsilon^\sigma,x_1+x_2\varepsilon,(1,\varepsilon^\sigma),(1,\varepsilon))\right)$$
$$+B_2(x_1)R_1[\varepsilon,\varepsilon^\sigma]+B_1(x_1)B_1(x_2)R_2[\varepsilon,\varepsilon^\sigma]+B_2(x_2)R_3[\varepsilon,\varepsilon^\sigma].$$
For $z=x_1+x_2 \varepsilon \in R(g)$, by the proof of Theorem 1 of \cite{ShintaniIII} the map
$$z \mapsto \overline{-z}:=\left\{
\begin{array}{ll}
1-x_1 & \textrm{ if $x_2=0$, $0 < x_1 < 1$}\\
1-x_1+(1-x_2)\varepsilon & \textrm{ if $0<x_1,x_2<1$}\\
1+(1-x_2)\varepsilon & \textrm{ if  $x_1=1$, $0 < x_2 < 1$}
\end{array} \right. $$
is a bijection between $R(g)$ and $R([\nu]_\mathfrak{f}^+g)$.
If $\overline {-z}=\overline x_1+\overline x_2 \varepsilon$, then observe that $B_2(x_1)=B_2(\overline x_1)$, $B_1(x_1)B_1(x_2)=B_1(\overline x_1)B_1(\overline x_2)$ and $B_2(x_2)=B_2(\overline x_2)$.  Hence when we compute $L_S^{(m)}(0,g)-L_S^{(m)}(0,[\nu]_\mathfrak{f}^+g)$ using the above expression, the presence of the terms with the polynomials $R_1,R_2$ and $R_3$ vanish.  We define
$$\mathbf{G}^r(z,\varepsilon):={G_2^r(z,(1,\varepsilon)) \over G_2^r(1+\varepsilon-z,(1,\varepsilon))},$$
$$
\mathbf{H}^{r,k}(z,v,\varepsilon,\varepsilon^\sigma):={H^{r,k}(z,v,(1,\varepsilon),(1,\varepsilon^\sigma)) \over H^{r,k}(1+\varepsilon-z,v,(1,\varepsilon),(1,\varepsilon^\sigma))}.
$$
It follows that $L_S^{(m)}(0,g)-L_S^{(m)}(0,[\nu]_\mathfrak{f}^+g)$ can be written as a finite linear combination with coefficients in $F$ of values of the form
$$
\begin{array}{c}
\log(N(\mathfrak{fa}_j))^{i}\mathbf{G}^r(x_1+x_2 \varepsilon,\varepsilon)\\
\log(N(\mathfrak{fa}_j))^{i}\mathbf{G}^r(x_1+x_2 \varepsilon^\sigma,\varepsilon^\sigma)\\
\log(N(\mathfrak{fa}_j))^{j}\mathbf{H}^{r,k}(x_1+x_2 \varepsilon,x_1+x_2 \varepsilon^\sigma,\varepsilon,\varepsilon^\sigma)\\
\log(N(\mathfrak{fa}_j))^{j}\mathbf{H}^{r,k}(x_1+x_2 \varepsilon^\sigma,x_1+x_2 \varepsilon,\varepsilon^\sigma,\varepsilon).
\end{array}$$
\end{proof}
The proof of Theorem \ref{special} in \S\ref{proof} would enable us to give an explicit expression for the value $L_S^{(m)}(\chi,0)$ as an element of of the field $K^{m}(F)$, which is transcendental over $F$.  By the proof of Theorem \ref{specialone} we could obtain an expression for 
$$L_S^{(m)}(0,g)-L_S^{(m)}(0,[\nu]_\mathfrak{f}^+g)$$ 
as an element of a field $L_\mathfrak{f}(F)$.  It would be interesting to investigate whether the field $L_\mathfrak{f}(F)$ is algebraic over $F$, or if not, whether any subfield of it was.  Indeed,  according to higher order versions of Stark's conjectures certain combinations (defined by the explicit expression for the L-value) of the special values of theta functions are strongly related to units in some algebraic extension of $F$.

\section{Appendix}

\subsection{The analyticity of $H^{q,k}$}\label{analH}
\begin{lemma}\label{conv}
Fix $v \in \B{C}$ and $\omega, \lambda \in \B{R}^2$ such that neither of the quotients $\omega_2/\omega_1$ and $\lambda_2/\lambda_1$ are negative, and assume that $\left|\omega_i \right|> \left|\lambda_i \right|$ for $i=1,2$.  Then there exists $R,r>1$ such that for $\left| t \right|$ sufficiently large and for all $N$
$$\left|C_N(t,v,\omega,\lambda)\right|\leq {1 \over r^N}{
\max\left\{e^{r\left|t \right|\left|\lambda_1+\lambda_2-v \right|} ,e^{-r\left|t \right|\left|\lambda_1+\lambda_2-v \right|} \right\} \over R^2}.$$
\end{lemma}
\begin{proof}
By the conditions on $\omega$ and $\lambda$, there exists $r>1$ such that the function $\mathfrak{g}(t,u,z,v,\omega,\lambda)$ defined in \eqref{g} is a meromorphic function in $u$ possessing no poles in the circle $\{\left| u\right|<r\} $ other than the one at zero.  By the definition of the functions $C_N(t,v,\omega,\lambda)$ in \eqref{coeffg}, by Cauchy's formula we have 
$$C_N(t,v,\omega,\lambda)={1 \over 2 \pi i} \oint_{\left|u \right|=r} {1 \over u^{N+1}}{ e^{(\left| \lambda \right|-v)tu} \over (1-e^{t(\omega_1+\lambda_1u)})(1-e^{t(\omega_2+\lambda_2u)})} du.$$
Therefore we obtain
$$\left| C_N(t,v,\omega,\lambda)\right|\leq {1 \over 2 \pi } {2 \pi r \over r^{N+1}}  \max_{\left|u \right|=r} \left\{{ e^{(\left| \lambda \right|-v)tu} \over (1-e^{t(\omega_1+\lambda_1u)})(1-e^{t(\omega_2+\lambda_2u)})}\right\}.$$
Now consider
$$\left|{ e^{(\left| \lambda \right|-v)tu} \over (1-e^{t(\omega_1+\lambda_1u)})(1-e^{t(\omega_2+\lambda_2u)})} \right|= 
{\left| e^{(\left| \lambda \right|-v)tu}\right| 
\over \left|1-e^{t(\omega_1+\lambda_1u)}\right|\left|1-e^{t(\omega_2+\lambda_2u)}\right|}
$$
$$\leq 
{\left| e^{(\left| \lambda \right|-v)tu}\right| 
\over \left|1-\left|e^{t(\omega_1+\lambda_1u)}\right|\right|\left|1-\left|e^{t(\omega_2+\lambda_2u)}\right|\right|} \qquad \qquad 
$$ 
\begin{equation} \label{qwerty}\qquad \qquad \leq 
{\max_{\left| u\right|=r}\left| e^{(\left| \lambda \right|-v)tu}\right| 
\over \min_{\left| u\right|=r}\left|1-\left|e^{t(\omega_1+\lambda_1u)}\right|\right| \min_{\left| u\right|=r}\left|1-\left|e^{t(\omega_2+\lambda_2u)}\right|\right|}.
\end{equation} 
The remainder of the proof is concerned with obtaining bounds for these maxima and minima.  On the circle we have $u=r^{i\theta}$ for $0 \leq \theta \leq 2 \pi$.  We first consider the denominator of \eqref{qwerty}, and put $t=t_1+it_2$:

$$\left|1-\left|e^{t(\omega_1+\lambda_1u)}\right|\right|=\left|1-e^{\Re(t(\omega_1+\lambda_1u))} \right|$$
$$=\left| 1-e^{t_1\omega_1+r t_1\lambda_1 \cos(\theta)-rt_2\lambda_1 \sin(\theta)} \right| $$  
This expression assumes its extremal values when $t_1 \cos(\theta)-t_2 \sin(\theta)$ does, which are equal to $\pm \left| t \right|$. 
Hence
$$\min_{\left|u \right|=r}\left|1-\left|e^{t(\omega_1+\lambda_1u)}\right|\right|=\min\left\{ \left|1-e^{t_1\omega_1+r \lambda_1\left|t\right|}\right|,\left|1-e^{t_1\omega_1-r \lambda_1\left|t\right|}\right|\right\},$$
and we obtain a similar expression for 
$$\min_{\left| u\right|=r}\left|1-\left|e^{t(\omega_2+\lambda_2u)}\right|\right|.$$
Therefore, there exists a $T \in\B{R}$ such that if $\left|t\right|>T$  then 
$$\min\left\{\min_{\left| u\right|=r}\left|1-\left|e^{t(\omega_1+\lambda_1u)}\right|\right| ,\min_{\left| u\right|=r}\left|1-\left|e^{t(\omega_2+\lambda_2u)}\right|\right|\right\} \geq R.$$
Now consider the numerator of \eqref{qwerty}, and put $v=v_1+iv_2$.  We calculate
$$
\begin{array}{rcl}
(\left|\lambda \right|-v)tu&=&(\lambda_1+\lambda_2-v_1-iv_2)(t_1+it_2)u\\
&=&\left(\left[(\lambda_1+\lambda_2-v_1)t_1+v_2t_2\right]+i\left[(\lambda_1+\lambda_2-v_1)t_2-t_1v_2 \right]\right)u
\end{array}
$$
The real part of the above expression is equal to 
$$r\left[(\lambda_1+\lambda_2-v_1)t_1+v_2t_2\right]\cos(\theta)-r\left[(\lambda_1+\lambda_2-v_1)t_2-t_1v_2 \right]\sin(\theta).$$
Hence the extremal values of 
$$\left|e^{(\left|\lambda \right|-v)tu}\right|=\left|e^{(\left|\lambda \right|-v)tr(\cos(\theta)+i \sin(\theta))}\right|$$
are $e^{\pm E}$ where 
$$
\begin{array}{rcl}
E^2&=&r^2\left[(\lambda_1+\lambda_2-v_1)t_1+v_2t_2\right]^2+r^2\left[(\lambda_1+\lambda_2-v_1)t_2-t_1v_2 \right]^2\\
&=&r^2\left((\lambda_1+\lambda_2-v_1)^2t_1^2 +v_2^2t_2^2+(\lambda_1+\lambda_2-v_1)^2t_2^2+t_1^1v_2^2\right)\\
&=&r^2\left((\lambda_1+\lambda_2-v_1)^2\left|t \right|^2 +v_2^2\left| t\right|^2\right)\\
&=&r^2 \left|t \right|^2 \left|\lambda_1+\lambda_2-v \right|^2
\end{array}
$$
\end{proof}

\begin{cor}\label{analext}
Under the same conditions as Lemma \ref{conv}, the integral
$$ \int_{I(\lambda,\infty)} e^{(\left|\omega \right|-z)t} J^k(\mathfrak{g}(u))(1) {\log(-t) \over t} dt$$
converges if $\Re(z)>r\left|\lambda_1+\lambda_2-v \right|+\omega_1+\omega_2$.
\end{cor}

\begin{proof}
The integrand is bounded on the circular path around the origin, so it suffices to show that the following integrals converge:
$$\int_\lambda^\infty e^{(\left|\omega \right|-z)t} J^k(\mathfrak{g}(u))(1) {\log(-t) \over t} dt;$$
$$\int_\lambda^\infty e^{(\left|\omega \right|-z)t} J^k(\mathfrak{g}(u))(1) {\log(-t)+2 \pi i \over t} dt.$$
Note that since $t$ dominates $\log(t)$ there exists $R_1$ such that if $t >R_1$.
$$\max\left\{\left|{\log(-t) \over t }\right|, \left|{\log(-t)+2 \pi i \over t} \right|\right\}<1.$$
Once again, on any finite interval $(\lambda,R)$ the integrands are bounded so it suffices to show that the following integral is convergent for sufficiently large $R$:

\begin{equation}\label{appen2}
\int_R^\infty \left|e^{(\left|\omega \right|-z)t} \right| \left|J^k(\mathfrak{g}(u))(1)\right| dt
\end{equation}
for sufficiently large $R$.
$$ \left|J^k(\mathfrak{g}(u))(1)\right| \leq \left|\sum_{N=1}^{\infty}{C_N(t,v,\omega,\lambda) \over N^k}\right|$$
$$ \leq \sum_{N=1}^{\infty}\left|{C_N(t,v,\omega,\lambda) }\right| $$
$$\leq  {\max\left\{e^{r\left|t \right|\left|\lambda_1+\lambda_2-v \right|} ,e^{-r\left|t \right|\left|\lambda_1+\lambda_2-v \right|} \right\} \over R^2}\sum_{N=1}^\infty{1 \over r^N}.$$
Then the expression of \eqref{appen2} is less than or equal to  
$${1 \over 1-r}{1 \over R^2}\int_R^\infty \left|e^{(\left|\omega \right|-z)t} \right| {\max\left\{e^{r\left|t \right|\left|\lambda_1+\lambda_2-v \right|} ,e^{-r\left|t \right|\left|\lambda_1+\lambda_2-v \right|} \right\} } dt $$
If $\Re(z)>r\left|\lambda_1+\lambda_2-v \right|+\omega_1+\omega_2$ then this integral converges.\\

\end{proof}

\begin{cor}
The integral \eqref{generalh} defines an analytic theta function for the pseudolattice $L=\B{Z}\omega_1+\B{Z}\omega_2$ in the region $\Re(z)>r\left|\lambda_1+\lambda_2-v \right|+\omega_1+\omega_2$.
\end{cor}
\begin{proof}
Using the integral formula we see that for $i=1,2$
$$\log(H^{k,q}(z+\omega_i,v,\omega,\lambda))$$
$$={1 \over 2 \pi i}\oint_{I(\lambda,\infty)} e^{(\left|\omega \right|-z)t} (1-e^{\omega_it})J^k(\mathfrak{g}(u))(1) {\log(-t) \over t} dt+\log(H^{k,q}(z,v,\omega,\lambda)).$$
In a similar way to the proof of Corollary \ref{analext} one can show that the above integral converges for $\Re(z)>r\left|\lambda_1+\lambda_2-v \right|+\omega_1+\omega_2$.
\end{proof}

\subsection{Contour Integrals used in \S\ref{proof}}

In this appendix we calculate two integrals that we used previously in \S\ref{proof} whilst proving Theorem \ref{special}.
\begin{lemma} \label{log1}
Let $m \in \B{Z}$,and $r \in \B{N}$. Then
$$\int_{I(\lambda,1)}u^m \log^r(u) du=\left\{ 
\begin{array}{ll}
\sum_{k=0}^{r-1} (-1)^k {r! \over (r-k)!} {(2 \pi i)^{r-k} \over (m+1)^{k+1}} & \textrm{if $m \neq -1$}\\
{(2 \pi i)^{r+1} \over r+1} & \textrm{if $m=-1$}.
\end{array}
\right.$$
\end{lemma}
\begin{proof}
Since the integrand is holomorphic at all points away from $0$, we know the value of the integral is independent of $\lambda$.  We will show that we can split this integral in to a finite sum of finite integrals, and take the limit as $\lambda \rightarrow 0$. We have 
$$\int_{I(\lambda,1)}u^m \log^r(u) du=\int_{1}^\lambda u^m \log^r(u)+\int_\lambda^1 u
^m \left[\log(u)+2 \pi i\right]^r$$
$$+i\lambda^{m+1}\int_{0}^{2 \pi} e^{(m+1)i\theta} \left[\log(\lambda)+2 \pi i \theta \right]^r d\theta.$$
The second integral is $O(\lambda)$.  The first integral is equal to 

\begin{equation}\label{binom1}
\sum_{k=0}^{r-1} \binom{r}{k} (2 \pi i)^{r-k}\int_\lambda^1 u^m \log^k(u) du.
\end{equation}
Let 
$$I_{m,k}=\int_\lambda^1 u^m \log(u)^k du$$
and note the decomposition $u^m \log(u)^k={u^{-1}} \times u^{m+1} \log(u)^k$.  Using integration by parts we obtain
$$I_{m,k}=\left[u^{m+1}\log(u)^{k+1} \right]_\lambda^1-(m+1)I_{m,k+1}-kI_{m,k}.$$
We are only interested in the limit of these integrals as $\lambda \rightarrow 0$.  Taking this limit we obtain
$$\lim_{\lambda \rightarrow 0} I_{m,k}=-{k \over m+1} \lim_{\lambda \rightarrow 0} I_{m,k-1}.$$
This yields 
$$\lim_{\lambda \rightarrow 0}I_{m,k}=(-1)^k {k! \over (m+1)^{k+1}}$$ which when substituted in to \eqref{binom1} yields the result.
\end{proof}

\begin{lemma}\label{polylog}
Let $d$ be a non-negative integer.  Then
$$\lim_{\lambda \rightarrow 0} \int_{\lambda}^1 {\log(u)^r \over a_1+a_2u} du={r! \over a_2}Li_{r+1}
\left(-{a_2 \over a_1} \right).$$
\end{lemma}
\begin{proof}
We first observe the identity

\begin{equation}\label{log}
\int_0^t {\log^r(u) \over a_1+a_2u}du={1 \over a_2} \sum_{i=0}^r (-1)^{i+1} Li_{i+1}(-at/b) \log(t)^{r-i} {r! \over (r-i)!}.
\end{equation}
This is easy to prove by differentiation, and using the identities
$$Li_1(z)=-\log(1-z) \qquad \textrm{and} \qquad Li_{s+1}(z)=\int_0^z {Li_s(t) \over t} dt.$$
Now split the integral in the statement of the lemma in to two parts:
$$\int_\lambda^1 {\log(u)^d \over a_1+a_2u} du=\int_0^1 {\log(u)^d \over a_1+a_2u} du-\int_0^\lambda {\log(u)^d \over a_1+a_2u} du.$$
Evaluating the integral of \eqref{log} at $t=1$ yields ${r! \over a_2}Li_{r+1}$.  The order of vanishing of $Li_j(t)$ is at least $1$ for $j \geq 1$. Hence $Li_j(\lambda) \log(\lambda)^n \rightarrow 0$ as $\lambda \rightarrow 0$
for any $j \geq 1, n \in \B{N}$.
\end{proof}


\end{document}